\documentclass[review]{elsarticle}
\usepackage{bm,amssymb,amsmath,mathrsfs}
\usepackage{amsthm}
\usepackage{lineno}
\usepackage[usenames]{color}

\usepackage{dcolumn}% Align table columns on decimal point

\def\half{{\textstyle\frac{1}{2}}}

\newcommand{\bq}{\begin{equation}}
\newcommand{\eq}{\end{equation}}
\newcommand{\bc}{\begin{center}}
\newcommand{\ec}{\end{center}}
\newcommand{\bit}{\begin{itemize}}
\newcommand{\eit}{\end{itemize}}
\newcommand{\ben}{\begin{enumerate}}
\newcommand{\een}{\end{enumerate}}

\theoremstyle{plain}
\newtheorem{theorem}{Theorem}[section]
\newtheorem*{theorem*}{Theorem}
\newtheorem{proposition}[theorem]{Proposition}

\newtheorem{remark}[theorem]{Remark}
\newtheorem{definition}[theorem]{Definition}
\newtheorem{conjecture}[theorem]{Conjecture}
\begin{document}

%\journal{\ }
%\journal{Expositiones Mathematicae}
\journal{(internal report CC23-6)}

\begin{frontmatter}

\title{Structure of the probability mass function of the Poisson distribution of order $k$}

\author[cc]{S.~R.~Mane}
\ead{srmane001@gmail.com}
\address[cc]{Convergent Computing Inc., P.~O.~Box 561, Shoreham, NY 11786, USA}

\begin{abstract}
The Poisson distribution of order $k$ is a special case of a compound Poisson distribution.
For $k=1$ it is the standard Poisson distribution.
Although its probability mass function (pmf) is known,
what is lacking is a {\em visual} interpretation, which a sum over terms with factorial denominators does not supply.
Unlike the standard Poisson distribution, the Poisson distribution of order $k$ can display a maximum of {\em four} peaks simultaneously,
as a function of two parameters: the order $k$ and the rate parameter $\lambda$.
This note characterizes the shape of the pmf of the Poisson distribution of order $k$.
The pmf can be partitioned into a single point at $n=0$, an increasing sequence for $n \in [1,k]$ and a mountain range for $n>k$ (explained in the text).
The ``parameter space'' of the pmf is mapped out and the significance of each domain is explained,
in particular the change in behavior of the pmf as a domain boundary is crossed.
A simple analogy (admittedly unrelated) is that of the discriminant of a quadratic with real coefficients:
its domains characterize the nature of the roots (real or complex),
and the domain boundary signifies the presence of a repeated root.
Something similar happens with the pmf of the Poisson distribution of order $k$.
As an application, this note explains the mode structure of the Poisson distribution of order $k$.
Improvements to various inequalities are also derived (sharper bounds, etc.).
New conjectured upper and lower bounds for the median and the mode are also proposed.
\end{abstract}

\vskip 0.25in


\end{frontmatter}

\newpage
\setcounter{equation}{0}
\section{\label{sec:intro} Introduction}
In two recent notes \cite{Mane_Poisson_k_CC23_3,Mane_Poisson_k_CC23_5},
the author presented numerical results for the Poisson distribution of order $k$ \cite{PhilippouGeorghiouPhilippou}.
It is a variant (or extension) of the well-known Poisson distribution.
We begin with its formal definition. 
\begin{definition}
  \label{def:pmf_Poisson_order_k}
  The Poisson distribution of order $k$ (where $k\ge1$ is an integer) and parameter $\lambda > 0$
  is an integer-valued statistical distribution with the probability mass function (pmf)
\bq
\label{eq:pmf_Poisson_order_k}
f_k(n;\lambda) = e^{-k\lambda}\sum_{n_1+2n_2+\dots+kn_k=n} \frac{\lambda^{n_1+\dots+n_k}}{n_1!\dots n_k!} \,, \qquad n=0,1,2\dots
\eq
\end{definition}
\noindent
For $k=1$ it is the standard Poisson distribution.
The Poisson distribution of order $k$ is a special case of the compound Poisson distribution introduced by Adelson \cite{Adelson1966}.
Although exact expressions for the mean and variance of the Poisson distribution of order $k$ are known \cite{PhilippouMeanVar},
exact results for its median and mode are difficult to obtain.

What is lacking is a {\em visual} interpretation of the pmf, which a formal mathematical sum 
such as eq.~\eqref{eq:pmf_Poisson_order_k} does not supply.
Unlike the standard Poisson distribution, the Poisson distribution of order $k$ can display a maximum of {\em four} peaks simultaneously,
as a function of two parameters: the order $k$ and the rate parameter $\lambda$.
This note characterizes the shape of the pmf of the Poisson distribution of order $k$.
The ``parameter space'' of the pmf is mapped out and the significance of each domain is explained,
in particular the change in behavior of the pmf as a domain boundary is crossed.
A simple analogy (admittedly unrelated) is that of the discriminant of a quadratic with real coefficients:
its domains characterize the nature of the roots (real or complex),
and the domain boundary signifies the presence of a repeated root.
Something similar happens with the pmf of the Poisson distribution of order $k$.
As an application, this note explains the mode structure of the Poisson distribution of order $k$.
Improvements to various inequalities are derived (sharper bounds, etc.).
New conjectured upper and lower bounds for the median and the mode are also proposed.

The structure of this paper is as follows.
Sec.~\ref{sec:notation} presents basic definitions and notation employed in this note.
Improvements to some published inequalities are also presented (sharper bounds, etc.)
Sec.~\ref{sec:pmf} quantifies the structure of the pmf of the Poisson distribution of order $k$.
Secs.~\ref{sec:median} and \ref{sec:mode} present new conjectures for upper and lower bounds for the median and mode, respectively.
Sec.~\ref{sec:conc} concludes.

\newpage
\setcounter{equation}{0}
\section{\label{sec:notation}Basic notation and definitions}
For later reference we define the parameter $\kappa=k(k+1)/2$.
We denote the mean by $\mu$, the median by $\nu$ and the mode by $m$ (with pertinent subscripts, etc.~to denote the dependence on $k$ and $\lambda$, see below).
Philippou \cite{PhilippouMeanVar} derived that the mean is $\mu_k(\lambda)=\kappa\lambda$ and the variance is $\sigma_k^2(\lambda) = \frac16 k(k+1)(2k+1)\lambda$.
For the median, we follow the exposition in \cite{AdellJodraPoisson2005}:
if $Y_{k,\lambda}$ is a random variable which is Poisson distributed with order $k$ and parameter $\lambda$,
the median is defined as the smallest integer $\nu$ such that $P(Y_{k,\lambda} \le \nu) \ge \frac12$.
With this definition, the median is unique and is always an integer.
The mode is defined as the location(s) of the {\em global maximum} of the probability mass function.
It is known that the mode may not be unique.
For the standard Poisson distribution with parameter $\lambda$, the mode equals $\lfloor\lambda\rfloor$ if $\lambda\not\in\mathbb{N}$,
but both $\lambda-1$ and $\lambda$ are modes if $\lambda\in\mathbb{N}$.
We adopt the following notation and definitions from \cite{KwonPhilippou}.
\begin{enumerate}
\item  
We work with $h_k(n;\lambda) = e^{k\lambda}f_k(n;\lambda)$ (see \cite{KwonPhilippou}) and refer to it as the ``scaled pmf'' below.
Observe from eq.~\eqref{eq:pmf_Poisson_order_k} that $h_k(n;\lambda)$ is a polynomial in $\lambda$ with all positive coefficients.
It has degree $n$ and for $n>0$ it has no constant term (also $h_k(0;\lambda)=1$ and $h_k(1;\lambda)=\lambda$ for all $k\ge1$).
Hence for fixed $k\ge1$ and $n>0$, $h_k(n;0)=0$ and $h_k(n;\lambda)$ is a strictly increasing function of $\lambda$ for $\lambda>0$.
\item
The parameter $r_k$ is defined as the positive root of the equation $h_k(k;\lambda)=1$.
It was shown in \cite{KwonPhilippou} that $r_k$ is unique and $0 < r_k < 1$.
\item
It was proved (Lemma 1 in \cite{KwonPhilippou}) that
for fixed $k\ge2$ and $\lambda>0$, the sequence $\{h_k(n;\lambda), n=1,\dots,k\}$ is strictly increasing,
i.e.~$h_k(n-1;\lambda) < h_k(n;\lambda)$ for $n=2,\dots,k$.
Hence only the last index $k$ can be a mode.
{\em An integer in the interval $[1,k-1]$ can never be a mode of the Poisson distribution of order $k$.}
\item
It was proved (Lemma 3 in \cite{KwonPhilippou}) that
for fixed $k\ge2$ and $0<\lambda\le r_k$, then $h_k(k;\lambda) > h_k(k+1;\lambda)$.
This makes $h_k(k;\lambda)$ a local maximum in the histogram of the pmf, for sufficiently small values of $\lambda$.
We shall see this below, when plotting graphs of the histogram of the pmf.
Note that the condition $0<\lambda\le r_k$ is sufficient {\em but not necessary} to attain $h_k(k;\lambda) > h_k(k+1;\lambda)$.
\end{enumerate}
By the term ``double mode'' we mean the distribution is bimodal, with joint modes at $m_1$ and $m_2$.
For $k=1$, the standard Poisson distribution, the integers $m_1$ and $m_2$ are always consecutive integers, but for $k\ge2$ this need not be so.
Kwon and Philippou \cite{KwonPhilippou} tabulated a list of double modes for $k=2,3,4$ and $0<\lambda\le2$.
It was shown in \cite{Mane_Poisson_k_CC23_3} that for any $k\ge2$,
the Poisson distribution of order $k$ has a denumerable infinity of double modes, consisting of pairs of consecutive integers.
A major topic of this note is to characterize the mode structure for $k\ge2$, including double modes with non-consecutive integers.
The existence of three or more joint modes is an open question.
This note presents numerical evidence that the Poisson distribution of order $k$ does not have three or more joint modes.

The term ``first double mode'' signifies the first time (smallest value of $\lambda$) that the Poisson distribution of order $k$ has a double mode.
The mode values in this case are $0$ and $m>0$.
The following notation was introduced in \cite{Mane_Poisson_k_CC23_5} for the first double mode.
The nonzero mode value was denoted by $\hat{m}_k$ and the corresponding value of $\lambda$ was denoted by $\hat\lambda_k$.
In terms of this notation, Kwon and Philippou \cite{KwonPhilippou} showed that $\hat{m}_k = k$ for $k=2,3,4$.
It was shown in \cite{Mane_Poisson_k_CC23_3} that $\hat{m}_k = k$ for $k=2,\dots,14$ and $\hat{m}_k > k$ for $k\ge15$.
It was proved in \cite{Mane_Poisson_k_CC23_5} that $\hat{m}_k \ge k$ for all $k\ge2$, but an exact formula for $k\ge15$ is not known to date.
Nor is it proved that $\hat{m}_k$ increases strictly with $k$ (for $k\ge15$)
although numerical tests up to at least $k=10^4$ have not found any exceptions \cite{Mane_Poisson_k_CC23_5}.

We take the opportunity here to present improvements to various published inequalities (sharper bounds, etc.).
The following inequality was derived in \cite{Mane_Poisson_k_CC23_5}
\bq
\label{eq:ineq_hatlam_rk}
1/\kappa \le \hat\lambda_k \le r_k < 1 \,.
\eq
Recall that it was proved in \cite{KwonPhilippou} that
for fixed $k\ge2$ and $\lambda>0$, the sequence $\{h_k(n;\lambda), n=1,\dots,k\}$ is strictly increasing.
It follows that the location of the first double mode is not less than $k$, i.e.~$\hat{m}_k\ge k$.
Next, the mode is bounded by the floor of the mean
(Theorem 2.1 in \cite{GeorghiouPhilippouSaghafi}, recall eq.~\eqref{eq:mode_Georghiou_etal_Thm2.1})
and Philippou \cite{PhilippouMeanVar} showed that the value of the mean is $\kappa\lambda$.
Hence
\bq
\label{eq:hatm_ineq_klamk}
\hat{m}_k \le \kappa\hat\lambda_k \,.
\eq
We can deduce two inequalities from this information.
\begin{proposition}
\label{prop:upbound_hatmk}
Using eq.~\eqref{eq:hatm_ineq_klamk} and $\hat\lambda_k < 1$, it follows that for all $k\ge2$,
\bq
k \le \hat{m}_k < \kappa \,.
\eq
\end{proposition}
\begin{proposition}
Using eq.~\eqref{eq:hatm_ineq_klamk} and $\hat{m}_k\ge k$, we deduce $\kappa\hat\lambda_k \ge k$.
Solving for $\hat\lambda_k$ yields the inequality
\bq
\label{prop:lowbound_hatlam}
\hat\lambda_k \ge \frac{k}{\kappa} = \frac{2}{k+1} \,.
\eq
\end{proposition}
\begin{remark}
Using eq.~\eqref{prop:lowbound_hatlam}, we improve the inequalities in eq.~\eqref{eq:ineq_hatlam_rk} as follows
\bq
\label{eq:ineq_hatlam_rk1}
\frac{2}{k+1} \le \hat\lambda_k \le r_k < 1 \,.
\eq
\end{remark}
\begin{remark}  
Philippou \cite{PhilippouFibQ} showed that the Poisson distribution of order $k$ has a unique mode of zero if $\lambda < 1/\kappa = 2/(k(k+1))$.
Using eq.~\eqref{prop:lowbound_hatlam}, we improve this bound to say the mode is uniquely zero if
\bq
\label{eq:newboundmodezero}
  \lambda < \frac{2}{k+1} \,.
\eq
This is a sufficient but not necessary condition.  
\end{remark}  
\noindent
\begin{remark}
We can now prove the conjecture in \cite{Mane_Poisson_k_CC23_3} that if the median is zero then the mode is also zero.
\end{remark}  
\begin{proof}
  The median is zero if and only if $\lambda \le (\ln2)/k$ (proved in \cite{Mane_Poisson_k_CC23_3}).
  Observe that for all $k\ge1$,
\bq
\frac{\ln2}{k} < \frac{2}{k+1} \,.
\eq
Hence if $\lambda \le (\ln2)/k$ (median is zero) then $\lambda < 2/(k+1)$ and from eq.~\eqref{eq:newboundmodezero} the mode is zero.
\end{proof}

\newpage
\setcounter{equation}{0}
\section{\label{sec:pmf}Structure of the probability mass function}
\subsection{\label{sec:pmf_gen}General remarks}
\begin{figure}[!htb]
\centering
\includegraphics[width=0.75\textwidth]{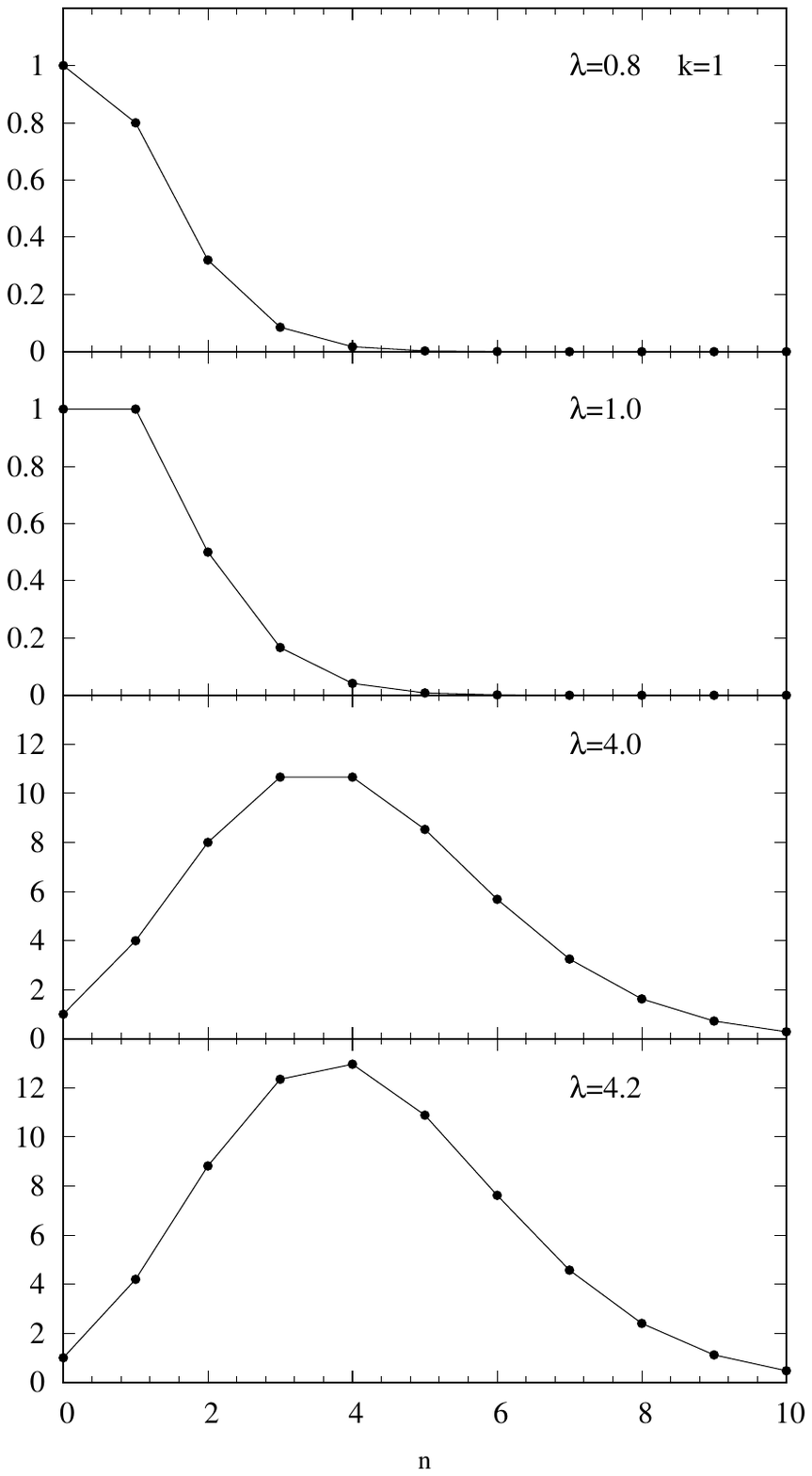}
\caption{\small
\label{fig:pmf_k1}
Plot of the scaled pmf $h_k(n;\lambda)$ of the standard Poisson distribution ($k=1$) for various values of the rate parameter $\lambda$.}
\end{figure}

\begin{figure}[!htb]
\centering
\includegraphics[width=0.75\textwidth]{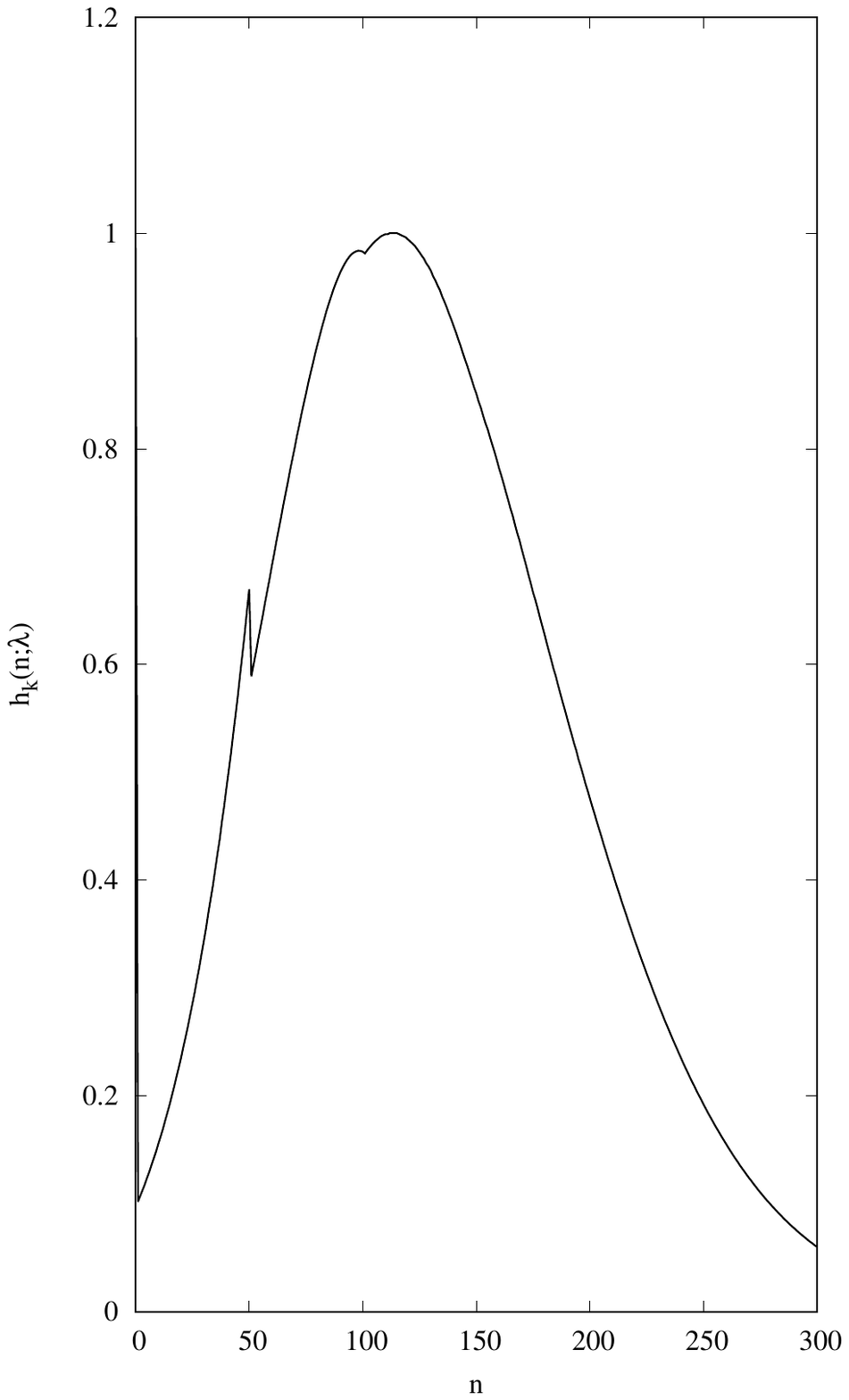}
\caption{\small
\label{fig:pmf_k50_lam_0_10194}
Plot of the scaled pmf $h_k(n;\lambda)$ for the Poisson distribution of order $50$ with rate parameter $\lambda=0.10194$.}
\end{figure}

We shall plot graphs to investigate the structure of the Poisson distribution of order $k$.
To fix ideas, we plot the value of the scaled pmf $h_k(n;\lambda)$, as opposed to the true pmf $f_k(n;\lambda)$ in eq.~\eqref{eq:pmf_Poisson_order_k}.
The prefactor $e^{-k\lambda}$ does not alter the shape of the histogram and is therefore irrelevant to the analysis below.
To demonstrate, the scaled pmf of the standard Poisson distribution is plotted in Fig.~\ref{fig:pmf_k1}.
\begin{enumerate}
\item
  For $\lambda < 1$, e.g.~$\lambda=0.8$, the unique mode is zero and the pmf decreases monotonically.
\item
  When $\lambda = 1$, the first double mode is attained, at $0$ and $1$.
\item
  As the value of $\lambda$ increases further, the pmf displays a single peak.
  Its location shifts to the right as the value of $\lambda$ increases.
  As the value of $\lambda$ increases, the location of the mode increases in unit steps.  
\item  
  If the value of $\lambda$ is an integer, e.g.~$\lambda=4$, the peak has a flat top and the distribution is bimodal, with joint modes at $\lambda-1$ and $\lambda$.
\item
  If the value of $\lambda$ is not an integer, e.g.~$\lambda=4.2$, the mode is unique, with value $\lfloor\lambda\rfloor$.
\end{enumerate}
By contrast, Fig.~\ref{fig:pmf_k50_lam_0_10194} displays a plot of the scaled pmf of the Poisson distribution of order $50$, for $\lambda\simeq0.10194$.
There are {\em four} peaks in the histogram,
at (i) $n=0$ (height$\;=1$),
(ii) $n=50$ (height$\;\simeq0.6698$), which is a local maximum and not a mode,
(iii) $n=98$ (height$\;\simeq0.98358$), which is also a local maximum and not a mode,
(iv) $n=113$ (height$\;=1$), which is a double mode along with $0$.
As the value of $\lambda$ changes, the relative heights of the peaks change, and the mode structure will vary.
For example, we know that for sufficiently small values $0 < \lambda < 2/(k+1)$, the unique mode is $0$.
Our goal in this section is to understand the evolution of the structure of the scaled pmf of the Poisson distribution of order $k$.

As complicated as it appears, Fig.~\ref{fig:pmf_k50_lam_0_10194} does contain useful clues as to how to parameterize the scaled pmf.
First, the point at $n=0$ always has a height of $1$ and is an invariant: it is the same for all $k$ and $\lambda$.
Next, it was mentioned earlier that the points for $1\le n \le k$ form a strictly increasing sequence 
(this was proved in Lemma 1 in \cite{KwonPhilippou}) and that sequence is visible in Fig.~\ref{fig:pmf_k50_lam_0_10194}.
It was also mentioned earlier that for sufficiently small $\lambda$, the point at $n=k$ is a local maximum 
(proved in Lemma 3 in \cite{KwonPhilippou}) and this fact is also visible in Fig.~\ref{fig:pmf_k50_lam_0_10194}.
The remaining points $n > k$ exhibit a more or less ``mountain range'' appearance, similar to the standard Poisson distribution,
except there are {\em two} peaks.
We shall refer to the region $n>k$ as the ``mountain range'' region and speak of the ``left peak'' and the ``right peak'' below.
However, be advised that the mountain range does not always have two peaks: that is why we need a ``parameter map'' to classify matters.

We wish to compose a ``parameter map'' to characterize the behavior of the scaled pmf of the Poisson distribution of order $k$
as a function of the order $k$ and the rate parameter $\lambda$.
For the dependence on $k$, there are four cases:
(i) $k \in [2,3]$, (ii) $k \in [4,14]$, (iii) $k \in [15,41]$, (iv) $k \ge 42$.
In each case, we increase the value of $\lambda$ continuously from zero and examine the structure of the scaled pmf, and the resulting consequences.

\newpage
\subsection{\label{sec:pmf_k2-k3}Case $k\in[2,3]$}
\begin{figure}[!htb]
\centering
\includegraphics[width=0.75\textwidth]{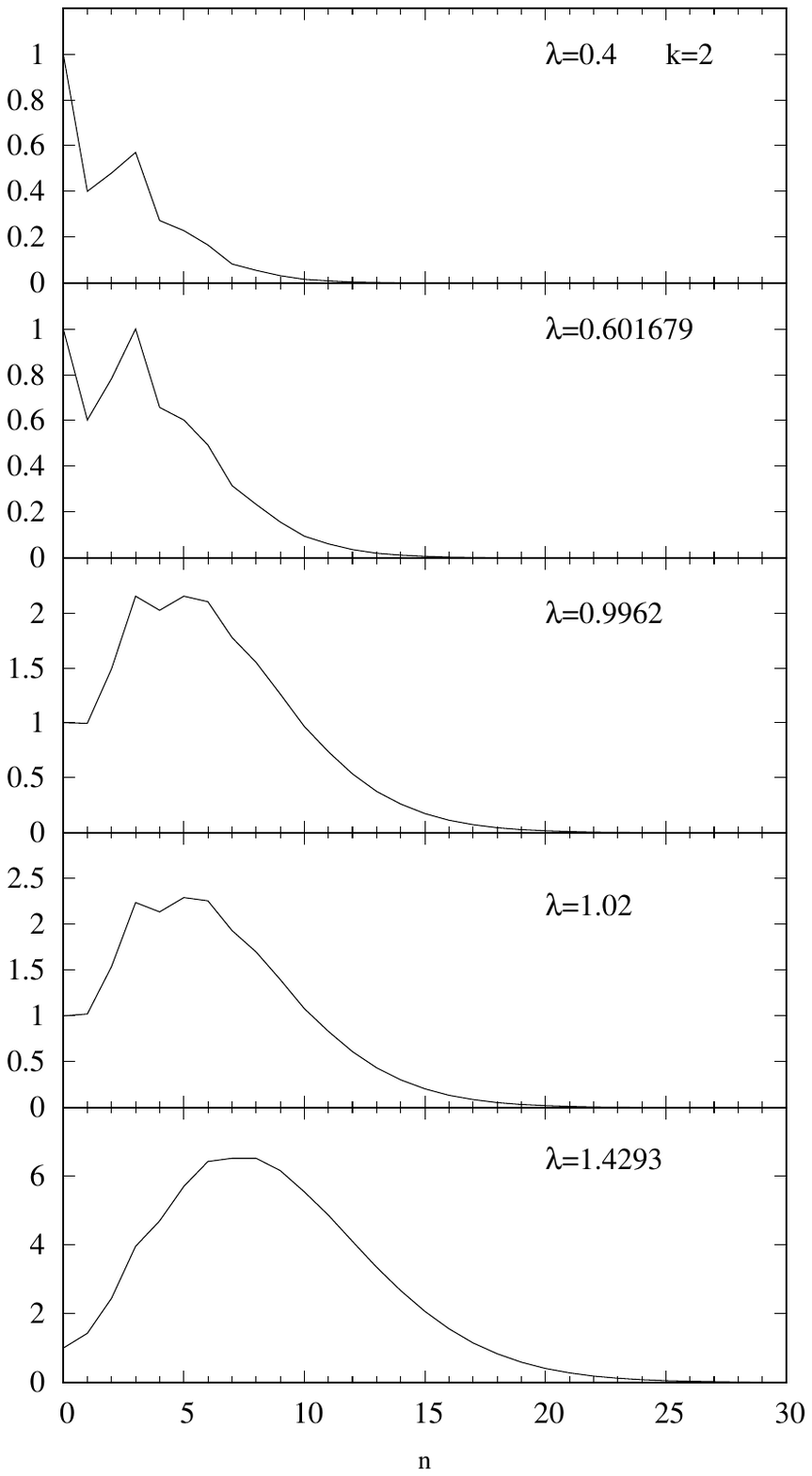}
\caption{\small
  \label{fig:pmf_k3}
Plot of the scaled pmf $h_k(n;\lambda)$ for the Poisson distribution of order $3$ for selected values of the rate parameter $\lambda$.}
\end{figure}

We select $k=3$ as a representative example, because for $k=2$ the ``increasing sequence'' $n=1,\dots,k$ contains only two points and is not informative.
The scaled pmf of the Poisson distribution of order $3$ is plotted in Fig.~\ref{fig:pmf_k3} for selected values of $\lambda$, increasing from top to bottom.
\begin{enumerate}
\item
  In the top panel, $\lambda=0.4$. This value is so small that the mountain range region ($n>k$) is monotonically decreasing.
  The increasing sequence for $n=1,\dots,k$ is visible, as is the local maximum at $n=k$, but it is not a global maximum.
  The unique mode is zero.
\item
  As the value of $\lambda$ increases, the height of the local maximum at $n=k$ reaches $1$ and a double mode is attained.
  The joint modes are at $0$ and $k$.
  This is displayed in the second panel, where $\lambda\simeq0.601679$.
  The mountain range region ($n>k$) is monotonically decreasing, although one can discern that a peak is forming.
\item
  When the value of $\lambda$ increases, the unique mode jumps from $0$ to $k$, i.e.~an increase of more than one unit.
  The ``first double mode'' represents a boundary between two domains in the parameter space.
\item
  As the value of $\lambda$ increases further, the mountain range develops a single peak,
  and its height rises to equal that of the point at $n=k$.
  This is displayed in the third panel, where $\lambda\simeq0.9962$.
  {\em The joint modes are at $k$ and $k+2$ (for both cases $k=2$ and $k=3$).}
\item
  When the value of $\lambda$ increases, the unique mode jumps from $k$ to $k+2$, i.e.~an increase of more than one unit.
  This is the ``second double mode'' and is a boundary between two domains in the parameter space.
\item
  As the value of $\lambda$ increases further, the height of the peak of the mountain range exceeds that of the point at $n=k$.
  The location of the single peak shifts rightwards and its height increases.
  The mode is determined solely by the location of the peak of the mountain range
  and increases in unit steps and takes all integer values $\ge k+2$.
  This is displayed in the fourth panel, where $\lambda=1.02$.
  The scaled pmf has a unique mode given by the height of the single mountain peak.
  The point at $n=k$ is still a local maximum but plays no further role to determine the mode.
\item  
  For a discrete (denumerably infinite) set of values of $\lambda$, the mountain peak has a flat top and the distribution is bimodal.
  The joint modes consist of consecutive integers.
  This is displayed in the fifth panel, where $\lambda\simeq1.4293$.
  There is a double mode, at $7$ and $8$.
  It is {\em almost} a triple mode, but the height of the point at $n=6$ is a little lower.
  The point at $n=k$ is no longer a local maximum.
\end{enumerate}

\newpage
\subsection{\label{sec:pmf_k4_14}Case $k\in[4,14]$}
\begin{figure}[!htb]
\centering
\includegraphics[width=0.75\textwidth]{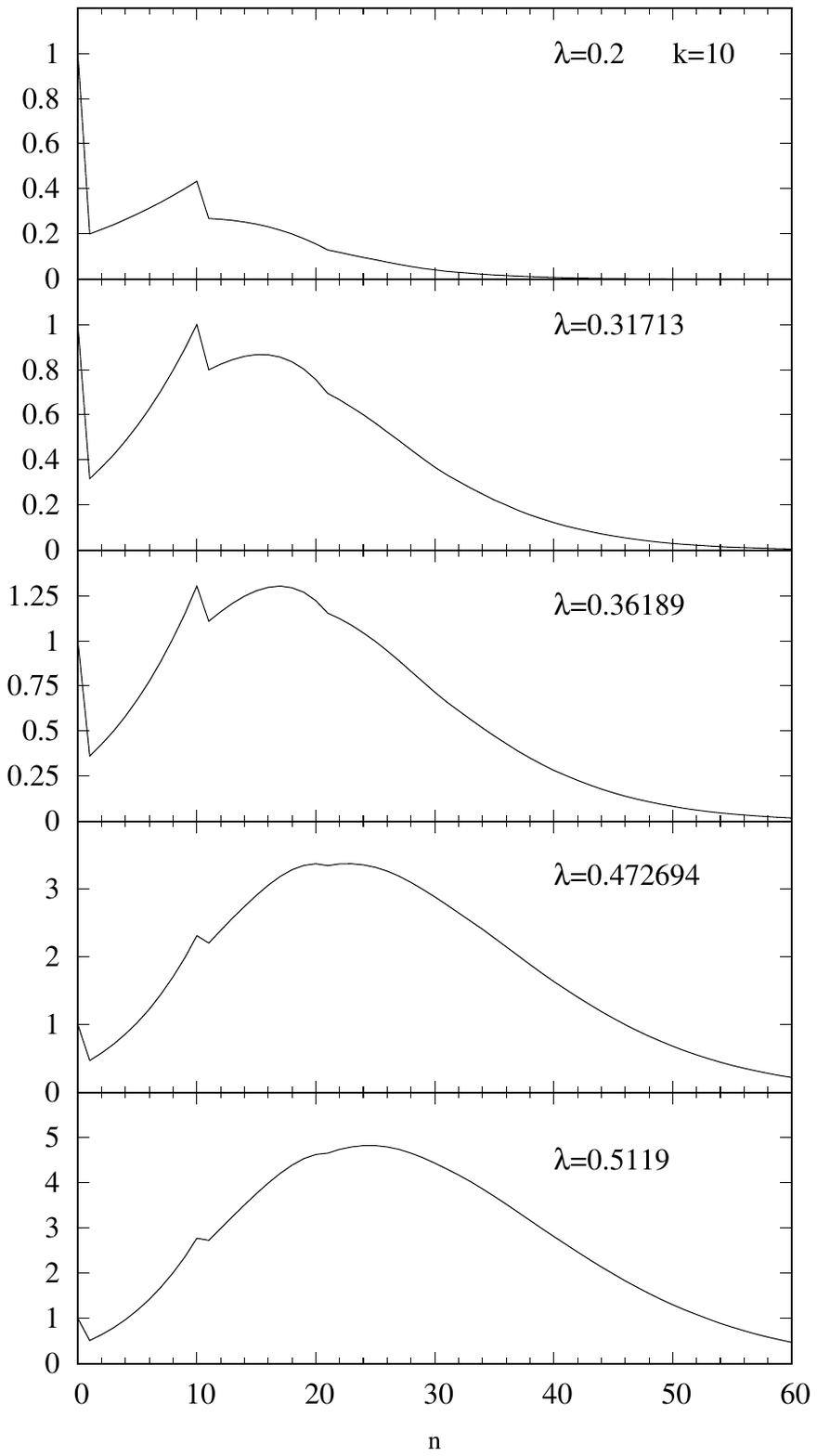}
\caption{\small
  \label{fig:pmf_k10}
Plot of the scaled pmf $h_k(n;\lambda)$ for the Poisson distribution of order $10$ for selected values of the rate parameter $\lambda$.}
\end{figure}

We select $k=10$ as a representative example.
The scaled pmf of the Poisson distribution of order $10$ is plotted in Fig.~\ref{fig:pmf_k10} for selected values of $\lambda$, increasing from top to bottom.
\begin{enumerate}
\item
  In the top panel, $\lambda=0.2$. This value is so small that the mountain range region ($n>k$) is monotonically decreasing.
  The increasing sequence for $n=1,\dots,k$ is visible, as is the local maximum at $n=k$, but it is not a global maximum.
  The unique mode is zero.
\item
  As the value of $\lambda$ increases, the height of the local maximum at $n=k$ reaches $1$ and a double mode is attained.
  The joint modes are at $0$ and $k$.
  This is displayed in the second panel, where $\lambda\simeq0.31713$.
  The mountain range region ($n>k$) exhibits a left peak, but its height is less than $1$ and is a local but not global maximum.
\item
  When the value of $\lambda$ increases, the unique mode jumps from $0$ to $k$, i.e.~an increase of more than one unit.
  The ``first double mode'' represents a boundary between two domains in the parameter space.
\item
  As the value of $\lambda$ increases further, the height of the left peak rises and equals that of the point at $n=k$.
  This is displayed in the third panel, where $\lambda\simeq0.36189$.
  The joint modes are at $k$ and the location of the left mountain peak, say $m_{\rm left}$.
  The value of $m_{\rm left}$ depends on $\lambda$ but is always at least $k+2$.
  No explicit formula is yet known for $m_{\rm left}$ when a double mode is attained with the point at $n=k$,
  although the values can be tabulated for all $k\in[4,14]$.
\item
  When the value of $\lambda$ increases, the unique mode jumps from $k$ to $m_{\rm left}$, i.e.~an increase of more than one unit.
  This is the ``second double mode'' and is a boundary between two domains in the parameter space.
\item
  As the value of $\lambda$ increases further, the mode is determined by the location of the left mountain peak.
  The value of the mode increases in unit steps {\em and there are double modes consisting of pairs of consecutive integers}.
  However, the mountain range develops a second (right) peak and its height rises faster than that of the left peak and it catches up with the height of the left peak.
\item
  As the value of $\lambda$ increases further, a double mode is attained, where the heights of the two mountain peaks are equal.
  This is displayed in the fourth panel, where $\lambda\simeq0.472694$.
  We denote the locations of the two mountain peaks by $m_{\rm left}$ and $m_{\rm right}$, respectively.
  {\em Very significant:} the value of $m_{\rm left}$ is {\em larger} than that in the previous panel, i.e.~the value of $m_{\rm left}$ increases with $\lambda$.
  The point at $n=k$ is still a local maximum but plays no further role to determine the mode.
  Observe also that there is a local minimum between the two mountain peaks, i.e.~they are separated by more than one unit and are not consecutive integers.
\item
  When the value of $\lambda$ increases, the unique mode jumps from $m_{\rm left}$ to $m_{\rm right}$, which is an increase of more than one unit.
  {\em However, this is not the ``third double mode'' because the value of $m_{\rm left}$ increased between the third and fourth panels,
    i.e.~there were double modes consisting of pairs of consecutive integers.}
  We may refer to this as the ``third mode jump'' (by which is meant an increase of the mode by more than one unit) and is a boundary between two domains in the parameter space.
  With this terminology, the first double mode is the {\em first mode jump} and the second double mode is the {\em second mode jump}.
\item
  As the value of $\lambda$ increases further, the height of the right mountain peak exceeds that of the left.
  The location of the right peak shifts rightwards and its height increases.
  The mode is determined solely by the location of the right mountain peak
  and increases in unit steps and takes all integer values from its value at the third mode jump upwards.
\item  
  For a discrete (denumerably infinite) set of values of $\lambda$, the right mountain peak has a flat top and the distribution is bimodal.
  The joint modes consist of consecutive integers.
  This is displayed in the fifth panel, where $\lambda\simeq0.5119$.
  There is a double mode at $24$ and $25$.
  The point at $n=k$ is still a local maximum but is no longer relevant to determine the mode.
\end{enumerate}

\newpage
\subsection{\label{sec:pmf_k15_41}Case $k\in[15,41]$}
\begin{figure}[!htb]
\centering
\includegraphics[width=0.75\textwidth]{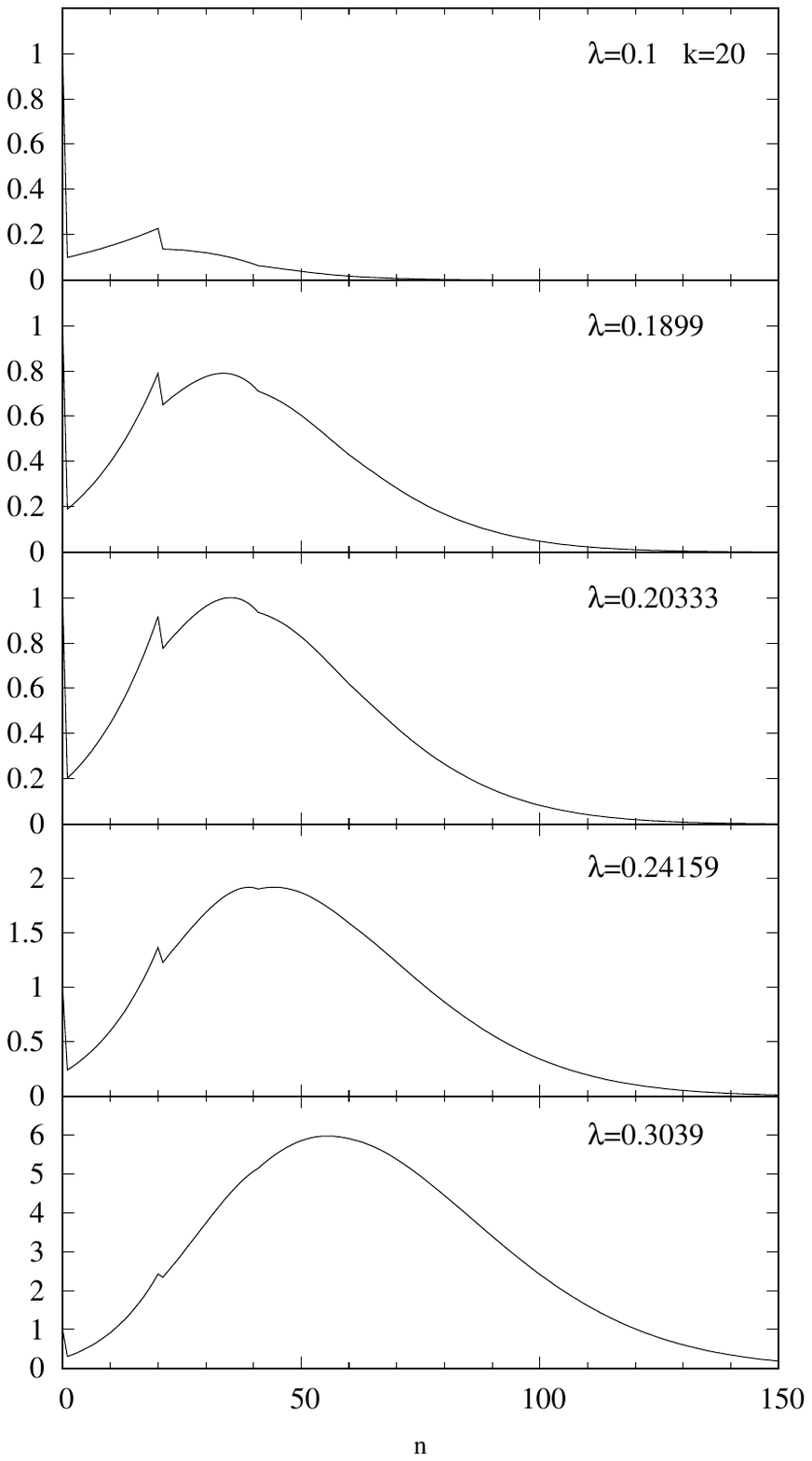}
\caption{\small
  \label{fig:pmf_k20}
Plot of the scaled pmf $h_k(n;\lambda)$ for the Poisson distribution of order $20$ for selected values of the rate parameter $\lambda$.}
\end{figure}

We select $k=20$ as a representative example.
The scaled pmf of the Poisson distribution of order $20$ is plotted in Fig.~\ref{fig:pmf_k20} for selected values of $\lambda$, increasing from top to bottom.
\begin{enumerate}
\item
  In the top panel, $\lambda=0.1$. This value is so small that the mountain range region ($n>k$) is monotonically decreasing.
  The increasing sequence for $n=1,\dots,k$ is visible, as is the local maximum at $n=k$, but it is not a global maximum.
  The unique mode is zero.
\item
  As the value of $\lambda$ increases, the mountain range develops a left peak and its height rises and catches up with that of the point at $n=k$.  
  This is displayed in the second panel, where $\lambda\simeq0.1899$.
  {\em Note that both heights are less than $1$, hence they are not modes.}
\item
  As the value of $\lambda$ increases further, the height of the left mountain peak rises faster than that of the point at $n=k$
  and it equals $1$ {\em before the point at $n=k$ does so.}
  This is displayed in the third panel, where $\lambda\simeq0.20333$.
  Hence the first double mode is {\em not} at $0$ and $k$.
  {\em The point at $n=k$ never plays a role to determine the mode.}
  The joint modes are at $0$ and the location of the left mountain peak, say $m_{\rm left}$.
  The value of $m_{\rm left}$ depends on $\lambda$ but is always at least $k+2$.
  No explicit formula is yet known for $m_{\rm left}$ when a double mode is attained with the point at $n=0$,
  although the values can be tabulated for all $k\in[15,41]$.
\item
  When the value of $\lambda$ increases, the unique mode jumps from $0$ to $m_{\rm left}$, i.e.~an increase of more than one unit.
  This is the ``first double mode'' (perhaps it is better to say the {\em first mode jump}) and is a boundary between two domains in the parameter space.
\item
  As the value of $\lambda$ increases further, the mode is determined by the location of the left mountain peak.
  The value of the mode increases in unit steps {\em and there are double modes consisting of pairs of consecutive integers}.
  However, the mountain range develops a second (right) peak and its height rises faster than that of the left peak and it catches up with the height of the left peak.
\item
  As the value of $\lambda$ increases further, a double mode is attained, where the heights of the two mountain peaks are equal.
  This is displayed in the fourth panel, where $\lambda\simeq0.24159$.
  We denote the locations of the two mountain peaks by $m_{\rm left}$ and $m_{\rm right}$, respectively.
  {\em Very significant:} the value of $m_{\rm left}$ is {\em larger} than that in the previous panel, i.e.~the value of $m_{\rm left}$ increases with $\lambda$.
  Observe also that there is a local minimum between the two mountain peaks, i.e.~they are separated by more than one unit and are not consecutive integers.
\item
  When the value of $\lambda$ increases, the unique mode jumps from $m_{\rm left}$ to $m_{\rm right}$, which is an increase of more than one unit.
  We refer to this as the ``second mode jump'' (by which is meant an increase of the mode by more than one unit) and is a boundary between two domains in the parameter space.
\item
  As the value of $\lambda$ increases further, the height of the right mountain peak exceeds that of the left.
  The location of the right peak shifts rightwards and its height increases.
  The mode is determined solely by the location of the right mountain peak
  and increases in unit steps and takes all integer values from its value at the second mode jump upwards.
\item  
  For a discrete (denumerably infinite) set of values of $\lambda$, the right mountain peak has a flat top and the distribution is bimodal.
  The joint modes consist of consecutive integers.
  This is displayed in the fifth panel, where $\lambda\simeq0.3039$.
  There is a double mode at $55$ and $56$.
\end{enumerate}

\newpage
\subsection{\label{sec:pmf_kge42}Case $k\ge42$}
\begin{figure}[!htb]
\centering
\includegraphics[width=0.75\textwidth]{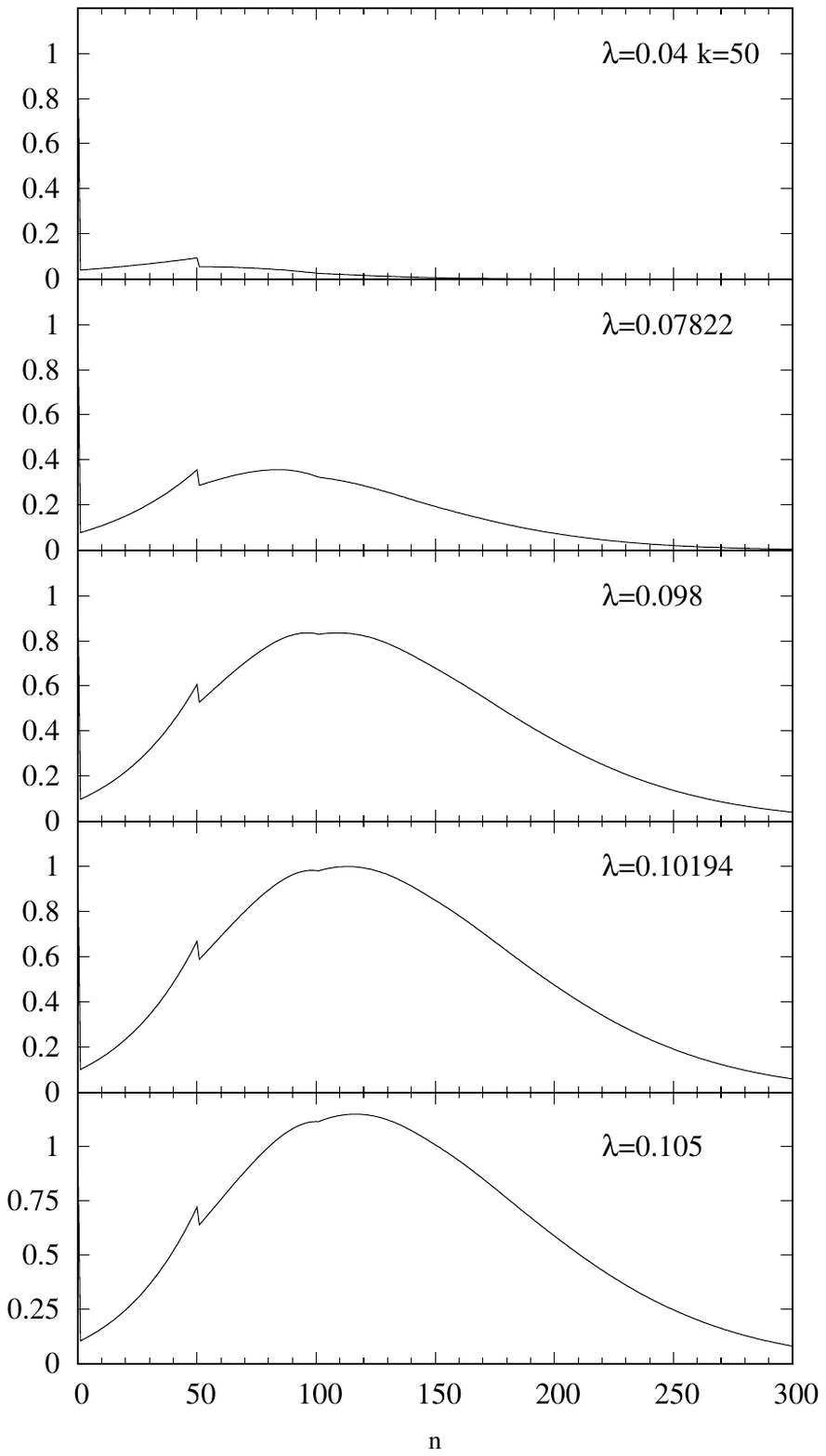}
\caption{\small
  \label{fig:pmf_k50}
Plot of the scaled pmf $h_k(n;\lambda)$ for the Poisson distribution of order $50$ for selected values of the rate parameter $\lambda$.}
\end{figure}

We select $k=50$ as a representative example.
The scaled pmf of the Poisson distribution of order $50$ is plotted in Fig.~\ref{fig:pmf_k50} for selected values of $\lambda$, increasing from top to bottom.
\begin{enumerate}
\item
  In the top panel, $\lambda=0.04$. This value is so small that the mountain range region ($n>k$) is monotonically decreasing.
  The increasing sequence for $n=1,\dots,k$ is visible, as is the local maximum at $n=k$, but it is not a global maximum.
  The unique mode is zero.
\item
  As the value of $\lambda$ increases, the mountain range develops a left peak and its height rises and catches up with that of the point at $n=k$.  
  This is displayed in the second panel, where $\lambda\simeq0.07822$.
  {\em Note that both heights are less than $1$, hence they are not modes.}
\item
  As the value of $\lambda$ increases further, the mountain range develops a right peak and its height rises and catches up with that of the left peak.
  This is displayed in the third panel, where $\lambda\simeq0.098$.
  {\em Note that both heights are less than $1$, hence they are not modes.}
\item
  As the value of $\lambda$ increases further, the height of the right mountain peak rises faster than that of the left peak
  and it equals $1$ {\em before the left peak or the point at $n=k$ does so.}
  This is displayed in the fourth panel, where $\lambda\simeq0.10194$.
  Hence the first double mode is {\em not} at $0$ and $k$.
  {\em The point at $n=k$ and the left mountain peak never play a role to determine the mode.}
  The joint modes are at $0$ and the location of the right mountain peak, say $m_{\rm right}$.
  The value of $m_{\rm right}$ depends on $\lambda$ but is always at least $k+2$.
  No explicit formula is yet known for $m_{\rm right}$ when a double mode is attained with the point at $n=0$,
  although the values can be tabulated for all desired values $k\ge42$.
\item
  When the value of $\lambda$ increases, the unique mode jumps from $0$ to $m_{\rm right}$, i.e.~an increase of more than one unit.
  This is the ``first double mode'' (perhaps it is better to say the {\em first mode jump}) and is a boundary between two domains in the parameter space.
\item
  As the value of $\lambda$ increases further, the mode is determined solely by the location of the right mountain peak.
  The location of the right peak shifts rightwards and its height increases.
  The mode increases in unit steps and takes all integer values from its value at the first mode jump upwards.
  {\em There is only one mode jump, for $k\ge42$.}
\item  
  For a discrete (denumerably infinite) set of values of $\lambda$, the right mountain peak has a flat top and the distribution is bimodal.
  The joint modes consist of consecutive integers.
  This is displayed in the fifth panel, where $\lambda\simeq0.105$.
  There is a double mode at $116$ and $117$.
\end{enumerate}

\newpage
\subsection{\label{sec:parametermap}Parameter map}
The shape of the pmf of the Poisson distribution of order $k$ can be partitioned into three sections:
(i) the single point at $n=0$, (ii) the increasing sequence ($n \in [1,k]$), and (iii) the mountain range ($n>k$).
It is by no means obvious from the formal sum in eq.~\eqref{eq:pmf_Poisson_order_k} that such a partition exists.
The increasing sequence is peculiar to the Poisson distribution of order $k$ and does not exist for the standard Poisson distribution ($k=1$).
It consists of terms where the number of summands equals $n$, because $n \in [1,k]$.
The mountain range consists of terms where the number of summands is strictly less than $n$, because $n>k$.
The mountain range is similar in character to the scaled pmf of the standard Poisson distribution.
For sufficiently small $\lambda>0$, it decreases monotonically, but for larger values of $\lambda$ it exhibits a peak.
It actually exhibits a maximum of {\em two} peaks.
The reason for this is not known: why not a single peak, and why not more than two peaks?

The mode structure of the scaled pmf of the Poisson distribution of order $k$ is determined by four parameters:
(i) the single point at $n=0$, whose height is always $1$, (ii) the single point at $n=k$, whose height depends on $\lambda$,
(iii) the left mountain peak $m_{\rm left}$, and
(iv) the right mountain peak $m_{\rm right}$.
For $k=2$ and $k=3$, there is only a single peak in the region $n>k$.
Both the location and the height of the mountain peaks depend on $\lambda$.
The evidence presented in this note suggests that the Poisson distribution of order $k$ does not have three or more joint modes.
The evidence is admittedly numerical, and cannot claim to be exhaustive, hence the existence of triple, etc.~modes remains open.

%\newpage
This is the parameter map:
\begin{enumerate}
\item
  For $k\in[2,3]$, the following happens as the value of $\lambda$ increases continuously from $0$.
  \begin{enumerate}
    \item
      For sufficiently small $\lambda>0$, the mode is zero.
    \item
      As $\lambda$ increases, the height of the point at $n=k$ reaches $1$ and there is a double mode at $0$ and $k$.
    \item
      When this domain boundary is crossed, the mode jumps from $0$ to $k$.
      This is the first mode jump.
    \item
      As $\lambda$ increases, the height of the single peak in the mountain range catches up to the height of the point at $n=k$ and there is a double mode at $k$ and $k+2$.
    \item
      When this domain boundary is crossed, the mode jumps from $k$ to $k+2$.
      This is the second mode jump.
    \item
      As $\lambda$ increases further, the mode is determined by the single peak in the mountain range.
      The location of the single peak shifts rightwards and its height increases.
      The mode increases in unit steps and takes all integer values $\ge k+2$.
      There is a denumerable infinity of double modes, consisting of pairs of consecutive integers.
  \end{enumerate}
\newpage
\item
  For $k\in[4,14]$, the following happens as the value of $\lambda$ increases continuously from $0$.
  \begin{enumerate}
    \item
      For sufficiently small $\lambda>0$, the mode is zero.
    \item
      As $\lambda$ increases, the height of the point at $n=k$ reaches $1$ and there is a double mode at $0$ and $k$.
    \item
      When this domain boundary is crossed, the mode jumps from $0$ to $k$.
      This is the first mode jump.
    \item
      As $\lambda$ increases, the height of the left mountain peak catches up to the height of the point at $n=k$ and there is a double mode at $k$ and $m_{\rm left}$.
      The value of $m_{\rm left}$ is not less than $k+2$.
    \item
      When this domain boundary is crossed, the mode jumps from $k$ to $m_{\rm left}$.
      This is the second mode jump.
    \item
      As $\lambda$ increases further, the mode is determined by the left mountain peak, which shifts rightwards and its height increases.
      The mode increases in unit steps and there are some double modes, consisting of pairs of consecutive integers.
    \item
      As $\lambda$ increases further, the height of the right mountain peak catches up to the height of the left mountain peak
      and there is a double mode at $m_{\rm left}$ and $m_{\rm right}$.
      The values of $m_{\rm left}$ and $m_{\rm right}$ are never consecutive integers.
    \item
      When this domain boundary is crossed, the mode jumps from $m_{\rm left}$ to $m_{\rm right}$.
      This is the third mode jump.
    \item
      As $\lambda$ increases further, the mode is determined by the right mountain peak, which shifts rightwards and its height increases.
      The mode increases in unit steps and takes all values from the third mode jump upwards.
      There is a denumerable infinity of double modes, consisting of pairs of consecutive integers.
  \end{enumerate}
\newpage
\item
  For $k\in[15,41]$, the following happens as the value of $\lambda$ increases continuously from $0$.
  \begin{enumerate}
    \item
      For sufficiently small $\lambda>0$, the mode is zero.
    \item
      As $\lambda$ increases, the height of the left mountain peak catches up to the height of the point at $n=k$,
      but their heights are less than $1$.
      The mode is zero.
    \item
      As $\lambda$ increases further, the height of the left mountain peak reaches $1$ and there is a double mode at $0$ and $m_{\rm left}$.
      The value of $m_{\rm left}$ is not less than $k+2$.
      The point at $n=k$ plays no role to determine the mode structure.
    \item
      When this domain boundary is crossed, the mode jumps from $0$ to $m_{\rm left}$.
      This is the first mode jump.      
    \item
      As $\lambda$ increases further, the mode is determined by the left mountain peak, which shifts rightwards and its height increases.
      The mode increases in unit steps and there are some double modes, consisting of pairs of consecutive integers.
    \item
      As $\lambda$ increases further, the height of the right mountain peak catches up to the height of the left mountain peak
      and there is a double mode $m_{\rm left}$ and $m_{\rm right}$.
      The values of $m_{\rm left}$ and $m_{\rm right}$ are never consecutive integers.
    \item
      When this domain boundary is crossed, the mode jumps from $m_{\rm left}$ and $m_{\rm right}$.
      This is the second mode jump.
    \item
      As $\lambda$ increases further, the mode is determined by the right mountain peak, which shifts rightwards and its height increases.
      The mode increases in unit steps and takes all values from the second mode jump upwards.
      There is a denumerable infinity of double modes, consisting of pairs of consecutive integers.
  \end{enumerate}
\newpage
\item
  For $k\ge42$, the following happens as the value of $\lambda$ increases continuously from $0$.
  \begin{enumerate}
    \item
      For sufficiently small $\lambda>0$, the mode is zero.
    \item
      As $\lambda$ increases, the height of the left mountain peak catches up to the height of the point at $n=k$,
      but their heights are less than $1$.
      The mode is zero.
    \item
      As $\lambda$ increases further, the height of the right mountain peak catches up to the height of the left mountain peak
      but their heights are less than $1$.
      The mode is zero.      
    \item
      As $\lambda$ increases further, the height of the right mountain peak reaches $1$ and there is a double mode at $0$ and $m_{\rm right}$.
      The value of $m_{\rm right}$ is not less than $k+2$.
      The point at $n=k$ and the left mountain peak play no role to determine the mode structure.
    \item
      When this domain boundary is crossed, the mode jumps from $0$ to $m_{\rm right}$.
      This is the first (and only) mode jump.
    \item
      As $\lambda$ increases further, the mode is determined by the right mountain peak, which shifts rightwards and its height increases.
      The mode increases in unit steps and takes all values from the first mode jump upwards.
      There is a denumerable infinity of double modes, consisting of pairs of consecutive integers.
  \end{enumerate}
\end{enumerate}

\newpage
\subsection{\label{sec:excluded}Excluded values}
We can now explain the ``excluded values'' tabulated in \cite{Mane_Poisson_k_CC23_5}, i.e.~integers which cannot be modes of the Poisson distribution of order $k$.
{\em The values tabulated in \cite{Mane_Poisson_k_CC23_5} are the integers which are skipped in the first, second or third mode jumps,}
and we can now explain the cause of those mode jumps.
A mode jump greater than unity occurs only when the controlling parameter of the mode changes,
e.g.~from $0$ to $k$ or from the left peak to the right peak, etc.
Only values $k \in [4,14]$ display three mode jumps, and for $k\ge42$ there is only one mode jump.
{\em But for all $k\ge2$ there is always a mode jump.}
As noted in \cite{Mane_Poisson_k_CC23_5}, the integers $1$ and $k+1$ are never modes of the Poisson distribution of order $k\ge2$.

\newpage
\setcounter{equation}{0}
\section{\label{sec:median}Median}
\begin{figure}[!htb]
\centering
\includegraphics[width=0.75\textwidth]{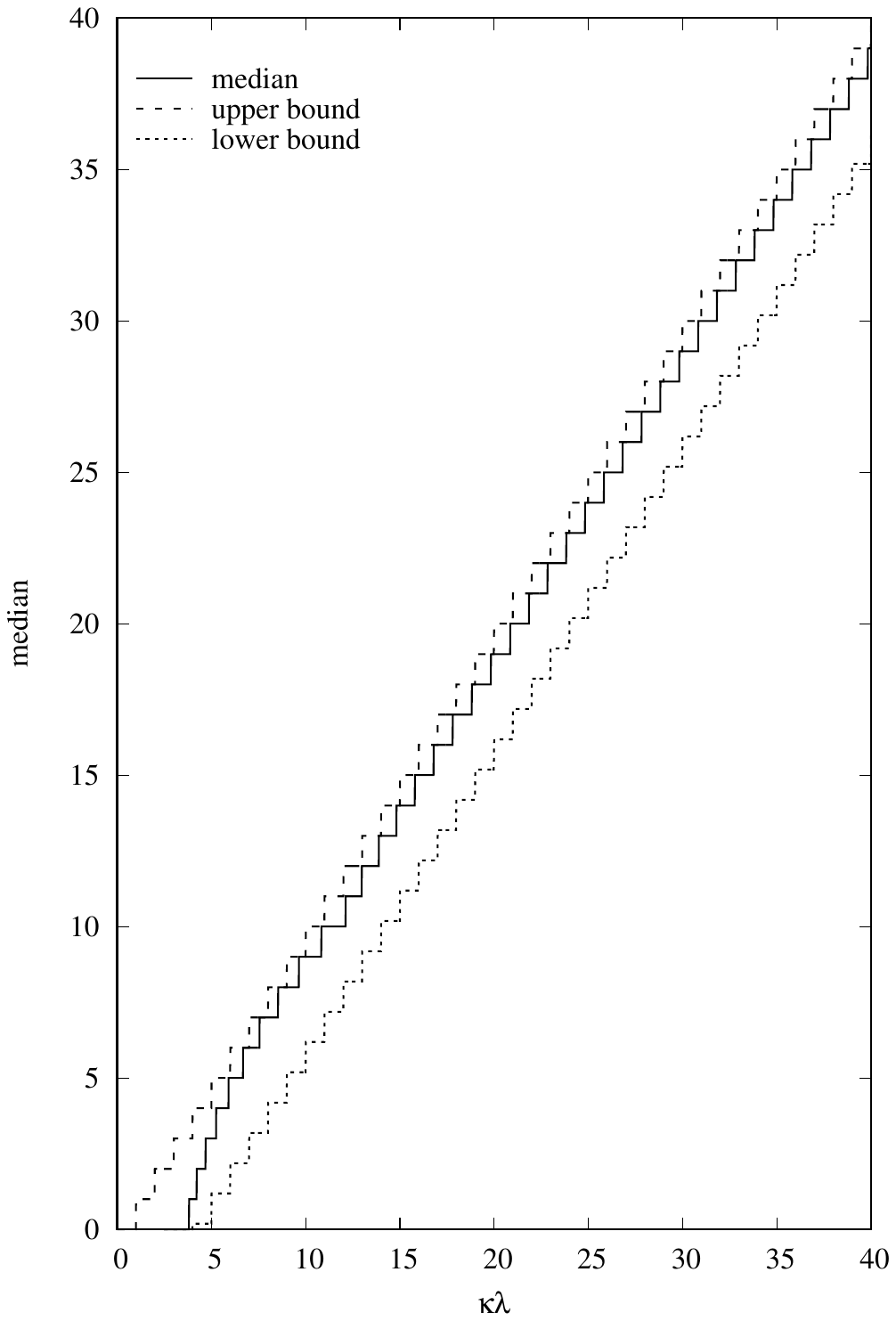}
\caption{\small
\label{fig:graph_median_k10}
Plot of the median (solid line) for the Poisson distribution of order $10$, with the upper bound (dashed) and lower bound (dotted).}
\end{figure}

It was proved in \cite{Mane_Poisson_k_CC23_3} that the median is zero if and only if $\lambda\le(\ln2)/k$.
Numerical studies reported in \cite{Mane_Poisson_k_CC23_3} also yielded the following expression for the median for $\lambda\ge1$.
If $n\in\mathbb{N}$ and $n\ge\kappa$, set $\lambda=n/\kappa$ (so $\lambda\ge1$).
Then the median is given by (eq.~(3.1) in \cite{Mane_Poisson_k_CC23_3})
\bq
\label{eq:Poisson_k_median_int_n}
\nu_k(n/\kappa) = n - \biggl\lfloor\frac{k+4}{8}\biggr\rfloor \,.
\eq
Note that $\frac{\kappa}{k}\,\ln2 = \frac12(k+1)\ln2 > \lfloor(k+4)/8\rfloor$ for all $k\ge1$.
Hence we conjecture the following bounds for the value of the median in the intermediate zone $\lambda\in(\frac1k\,\ln2,\,1)$.
\begin{conjecture}
  For fixed $k\ge2$ and $\lambda\in(\frac1k\,\ln2,\,1)$, we claim the lower bound for the median is 
\bq
\label{eq:non_asymp_median_lower_bound}
\nu_k(\lambda) \ge \max\bigl\{0, \lfloor\kappa\lambda\rfloor - \half (k+1)\ln2 \bigr\} \,.
\eq
For an upper bound, we know from the numerical studies reported in \cite{Mane_Poisson_k_CC23_3}
it is possible for the value of the median to exceed the mean, because $\lfloor(k+4)/8\rfloor=0$ if $k < 4$.
We propose the following upper bound for the median, where $c_0=1$ if $k \le3$ and $c_0=0$ if $k\ge4$.
\bq
\label{eq:non_asymp_median_upper_bound}
\nu_k(\lambda) \le \lfloor\kappa\lambda\rfloor + c_0 \,.
\eq
\end{conjecture}
\noindent
Monte Carlo scans using values $2\le k \le 4\times10^4$ and $\lambda \in (\frac1k\,\ln2,\,1)$ found
no violations of the bounds in eqs.~\eqref{eq:non_asymp_median_lower_bound} or \eqref{eq:non_asymp_median_upper_bound}.
Note the following:
\begin{enumerate}
\item
  From eq.~\eqref{eq:Poisson_k_median_int_n}, it is tempting to conjecture that a sharper upper bound is
  $\nu_k(\lambda) \le \lfloor\kappa\lambda\rfloor - \lfloor(k+4)/8\rfloor$
  or possibly $\nu_k(\lambda) \le 1 + \lfloor\kappa\lambda\rfloor - \lfloor(k+4)/8\rfloor$,
  but they are both false for $\lambda<1$.
  An example is $k=514$ and $\lambda\simeq0.0031619$,
  whence $\kappa\lambda \simeq 418.4998$ and $\lfloor\kappa\lambda\rfloor - \lfloor(k+4)/8\rfloor = 418 - 64 = 354$,
  but the median value is $367$.
\item
  A graph of the median is plotted in Fig.~\ref{fig:graph_median_k10}, for $k=10$ and $0 < \kappa\lambda \le 40$,
  whence $\kappa=55$ and $0 < \lambda \le 40/55 = 0.7272\dots$
  The value of the median is plotted as the solid line.
  The lower bound from eq.~\eqref{eq:non_asymp_median_lower_bound} is plotted as the dotted line.
  The upper bound from eq.~\eqref{eq:non_asymp_median_upper_bound} is plotted as the dashed line.
  Observe that the upper bound equals the median at several places, indicating that it is a sharp upper bound.
  Observe also that the lower bound is essentially a straight line,
  and does not capture the curvature as the value of the median approaches zero for small values of $\kappa\lambda$.
  This indicates that the expression for the lower bound can be improved.
\item
  The graph of the median in Fig.~\ref{fig:graph_median_k10} looks hyperbolic, with a straight line asymptote for large values of $\kappa\lambda$
  (see eq.~\eqref{eq:Poisson_k_median_int_n} for $\lambda \ge1$).
  This suggests a better parameterization for the value of the median is possible in the interval $\lambda\in(\frac1k\ln2,\,1)$.
  The matter is left for future work.
\end{enumerate}
Given the above, we propose the following exact results and bounds for the median, for all $k\ge2$ and $\lambda>0$.
In all cases, we fix $k\ge2$ and the median is denoted by $\nu_k(\lambda)$.
\begin{enumerate}
\item
  For $\lambda \in (0,\frac1k\,\ln2]$ the median is zero: $\nu_k(\lambda) = 0$. This was proved in \cite{Mane_Poisson_k_CC23_3}.
\item
For $\lambda \in (\frac1k\,\ln2,\,1)$, the median is bounded via eqs.~\eqref{eq:non_asymp_median_lower_bound} and \eqref{eq:non_asymp_median_upper_bound},
where $c_0=1$ if $k \le3$ and $c_0=0$ if $k\ge4$.
\bq
\max\bigl\{0, \lfloor\kappa\lambda\rfloor - \half (k+1)\ln2 \bigr\} \le
\nu_k(\lambda) \le \lfloor\kappa\lambda\rfloor + c_0 \,.
\eq
\item
For $\lambda \ge 1$, the median is given by eqs.~(3.1), (3.3) and (3.4) in \cite{Mane_Poisson_k_CC23_3} as follows.
Let $n\in\mathbb{N}$ and $n\ge\kappa$, then for $\kappa\lambda \in (\alpha_{k,n-1},\alpha_{k,n}]$ the median equals $\nu_k(n/\kappa)$ as follows.
\begin{subequations}
\begin{align}
\nu_k(n/\kappa) &= n - \biggl\lfloor\frac{k+4}{8}\biggr\rfloor \,,
\\
\alpha_{k,n} &= n + \textrm{frac}\Bigl(\frac{k+4}{8}\Bigr) +\frac{k}{8(2k+1)} + A_{k,n} \,,
\\
A_{k,n} &= \biggl(\frac{3\kappa}{349} + \frac{13}{1000}\biggr)\frac{1}{n} + \frac{13}{1500}\biggl(\biggl\lfloor\frac{k+4}{8}\biggr\rfloor -3\biggr)\,\frac{\kappa}{n^2} +\cdots
\end{align}
\end{subequations}
The expression for $A_{k,n}$ is approximate,
but the expressions for $\nu_k(n/\kappa)$ and $\alpha_{k,n}$ are otherwise conjectured to be exact.
\end{enumerate}

\newpage
\setcounter{equation}{0}
\section{\label{sec:mode}Mode}
Recall $\hat{m}_k$ is the value of the first double mode of the Poisson distribution of order $k$,
i.e.~the distribution is bimodal, with modes at $0$ and $\hat{m}_k$
and $\hat\lambda_k$ is the corresponding value of $\lambda$.
It was shown in \cite{Mane_Poisson_k_CC23_5} that $\hat\lambda_k$ is a strictly decreasing function of $k$ and also that $\hat{m}_k \ge k$.

\begin{figure}[!htb]
\centering
\includegraphics[width=0.75\textwidth]{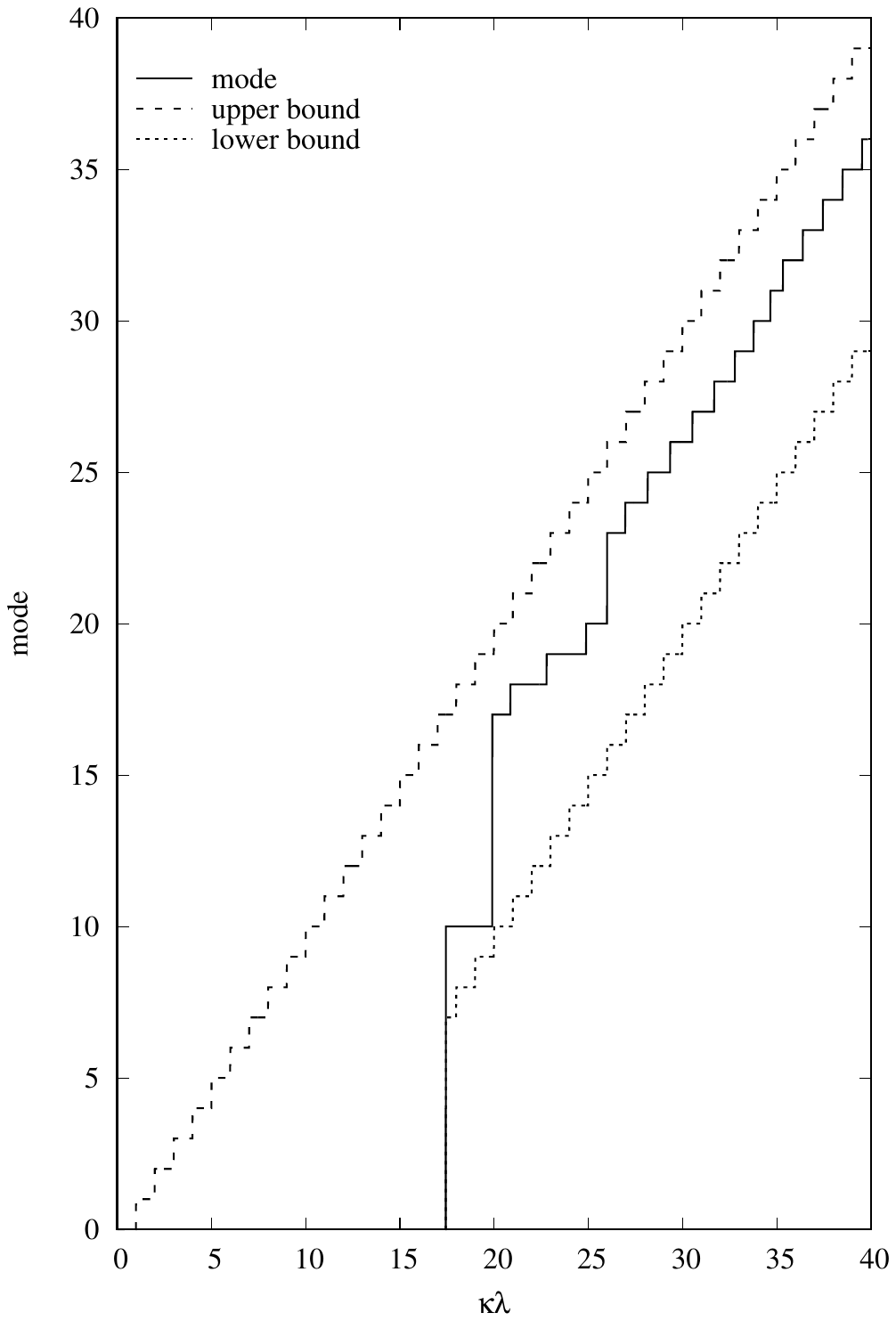}
\caption{\small
\label{fig:graph_mode_k10}
Plot of the mode (solid line) for the Poisson distribution of order $10$, with the upper bound (dashed) and lower bound (dotted).}
\end{figure}

\begin{figure}[!htb]
\centering
\includegraphics[width=0.75\textwidth]{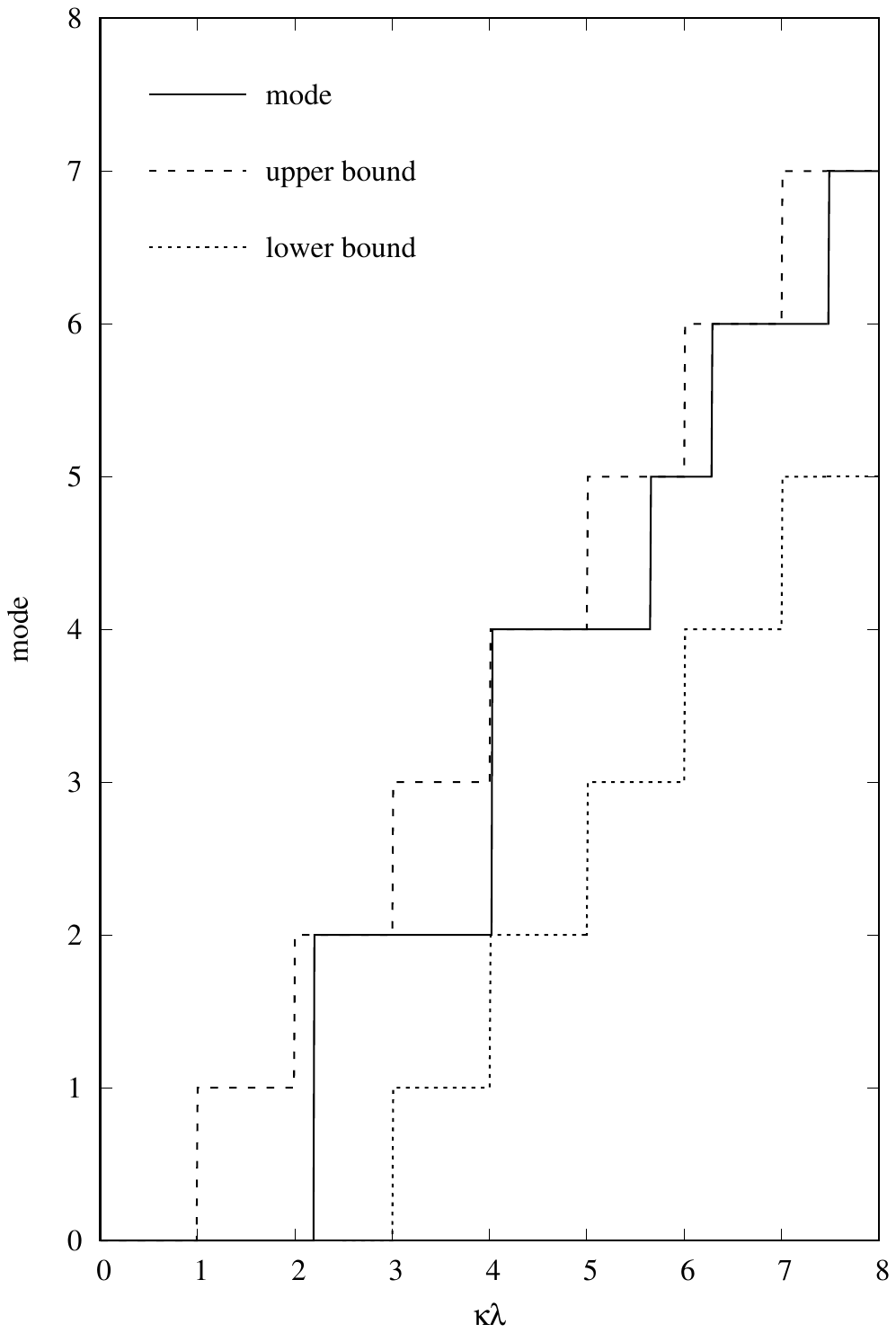}
\caption{\small
\label{fig:graph_mode_k2}
Plot of the mode (solid line) for the Poisson distribution of order $2$, with the upper bound (dashed) and lower bound (dotted).}
\end{figure}

The following asymptotic expression for the mode was conjectured in \cite{Mane_Poisson_k_CC23_3}.
Let $n\in\mathbb{N}$ and $n\ge2\kappa$ and set $\lambda=n/\kappa$ (so $\lambda\ge2$).
Then the mode is given by (eq.~(4.1) in \cite{Mane_Poisson_k_CC23_3})
\bq
\label{eq:Poisson_k_mode_int_n}
m_k(n/\kappa) = n - \biggl\lfloor\frac{3k+5}{8}\biggr\rfloor \,.
\eq
Theorem 2.1 in \cite{GeorghiouPhilippouSaghafi} states the following upper and lower bounds for the mode.
\bq  
\label{eq:mode_Georghiou_etal_Thm2.1}
\bigl\lfloor \kappa\lambda \bigr\rfloor - \kappa + 1 - \delta_{k,1} \le m_k(\lambda) \le \bigl\lfloor \kappa\lambda \bigr\rfloor \,.
\eq
We conjecture an improved lower bound for the mode.
\begin{conjecture}
  \label{conj:mode_low_bound}
  For fixed $k\ge2$ and $\lambda\in(\hat\lambda_k,\,2)$ (so the value of the mode is nonzero),
  we propose the following as an improved lower bound for the mode
\bq
\label{eq:non_asymp_mode_lower_bound}
m_k(\lambda) \ge \lfloor\kappa\lambda\rfloor -k \,.
\eq
The right-hand side is nonnegative and evidence will be presented below that it is a sharp lower bound.
\end{conjecture}
\noindent
Monte Carlo scans using values $2\le k \le 4\times10^4$ and $0 < \lambda < 2$ found no violations of eq.~\eqref{eq:non_asymp_mode_lower_bound}.
(Only cases where the value of the mode was positive were included in the scan, to satisfy the requirements of Conjecture \ref{conj:mode_low_bound}.)
Note the following:
\begin{enumerate}
\item
  From eq.~\eqref{eq:Poisson_k_mode_int_n}, it is tempting to conjecture that a sharper upper bound is
  $m_k(\lambda) \le \lfloor\kappa\lambda\rfloor - \lfloor(3k+5)/8\rfloor$
  or possibly $m_k(\lambda) \le 1 + \lfloor\kappa\lambda\rfloor - \lfloor(3k+5)/8\rfloor$,
  but they are both false for $\lambda<2$.
  An example is $k=44$ and $\lambda\simeq0.114198$,
  whence $\kappa\lambda \simeq 113.056$ and $\lfloor\kappa\lambda\rfloor - \lfloor(3k+5)/8\rfloor = 113 - 17 = 96$,
  but the mode value is $98$.
\item
  Note however that Monte Carlo scans have thus far failed to find an example where
  $m_k(\lambda) > 0$ and $m_k(\lambda) \ge 3 + \lfloor\kappa\lambda\rfloor - \lfloor(3k+5)/8\rfloor$.
\item
  A graph of the mode is plotted in Fig.~\ref{fig:graph_mode_k10}, for $k=10$ and $0 < \kappa\lambda \le 40$,
  whence $\kappa=55$ and $0 < \lambda \le 40/55 = 0.7272\dots$
  The value of the mode is plotted as the solid line.
  The lower bound from eq.~\eqref{eq:non_asymp_mode_lower_bound} is plotted as the dotted line.
  The upper bound is $\lfloor\kappa\lambda\rfloor$ (from \cite{GeorghiouPhilippouSaghafi}) and is plotted as the dashed line.
\item
  Unlike the median, the value of the mode does not always increase in unit steps.
  Observe in Fig.~\ref{fig:graph_median_k10} that the value of the mode jumps from $0$ to $10$ and then from $10$ to $17$ and then from $20$ to $23$.
  These mode jumps were explained in Sec.~\ref{sec:pmf}.
\item
  In Fig.~\ref{fig:graph_median_k10}, it looks as if the lower bound equals the mode at $\kappa\lambda = 20$, but it does not.
  The value of the mode jumps from $10$ to $17$ at $\kappa\lambda \simeq 19.91$,
  while the lower bound increases from $9$ to $10$ at $\kappa\lambda = 20$.
\item
  However, the case $k=2$ demonstrates that the expression for the lower bound in eq.~\eqref{eq:non_asymp_mode_lower_bound} is a sharp lower bound.
  A graph of the mode is plotted in Fig.~\ref{fig:graph_mode_k2}, for $k=2$ and $0 < \kappa\lambda \le 8$,
  whence $\kappa=3$ and $0 < \lambda \le 8/3 = 2.666\dots$
  The value of the mode is plotted as the solid line.
  The upper and lower bounds are plotted as the dashed and dotted lines, respectively.
  Both the mode and lower bound equal $2$ at $\kappa\lambda=4$.
  For $\kappa\lambda = 4.0238$ the value of the mode jumps to $4$ but the lower bound remains at $2$.
  Hence the mode and lower bound both equal $2$ in a nonempty subset of the interval $\kappa\lambda \in [4,4.0238)$.
\end{enumerate}
Given the above, we propose the following exact results and bounds for the mode, for all $k\ge2$ and $\lambda>0$.
In all cases, we fix $k\ge2$ and the mode is denoted by $m_k(\lambda)$.
\begin{enumerate}
\item
  For $\lambda \in (0,\hat\lambda_k)$ the mode is zero: $m_k(\lambda) = 0$, by the definition of $\hat\lambda_k$.
\item
  For $\lambda=\hat\lambda_k$, there is a double mode at $0$ and $\hat{m}_k$ (this is the definition of $\hat{m}_k$).
\item
For $\lambda \in(\hat\lambda_k,\,2)$, the mode is bounded as follows
\bq
\lfloor\kappa\lambda\rfloor -k \le m_k(\lambda) \le \lfloor\kappa\lambda\rfloor \,.
\eq
The lower bound is based on numerical studies in this note.
The upper bound is from \cite{GeorghiouPhilippouSaghafi}.
\item
For $\lambda \ge 2$, the mode is given by eqs.~(4.1), (4.5) and (4.6) in \cite{Mane_Poisson_k_CC23_3} as follows.
Let $n\in\mathbb{N}$ and $n\ge2\kappa$, then for $\kappa\lambda \in (\beta_{k,n-1},\beta_{k,n}]$ the mode equals $m_k(n/\kappa)$ as follows.
\begin{subequations}
\begin{align}
m_k(n/\kappa) &= n - \biggl\lfloor\frac{3k+5}{8}\biggr\rfloor \,,
\\
\beta_{k,n} &= n + \textrm{frac}\Bigl(\frac{3k+5}{8}\Bigr) +\frac{k-1}{8(2k+1)} + B_{k,n} \,,
\\
B_{k,n} &= \biggl(\frac{\kappa}{16 +\frac89} - \frac{1}{13+\frac23}\biggr)\frac{1}{n}
+ \biggl\lfloor\frac{3k+5}{8}\biggr\rfloor \,\frac{3\kappa}{50n^2} +\cdots
\end{align}
\end{subequations}
The expression for $B_{k,n}$ is approximate,
but the expressions for $m_k(n/\kappa)$ and $\beta_{k,n}$ are otherwise conjectured to be exact.
\item
Note that no general formula is yet known for $\hat{m}_k$ or $\hat\lambda_k$, although specific cases have been solved.
\end{enumerate}

\begin{figure}[!htb]
\centering
\includegraphics[width=0.75\textwidth]{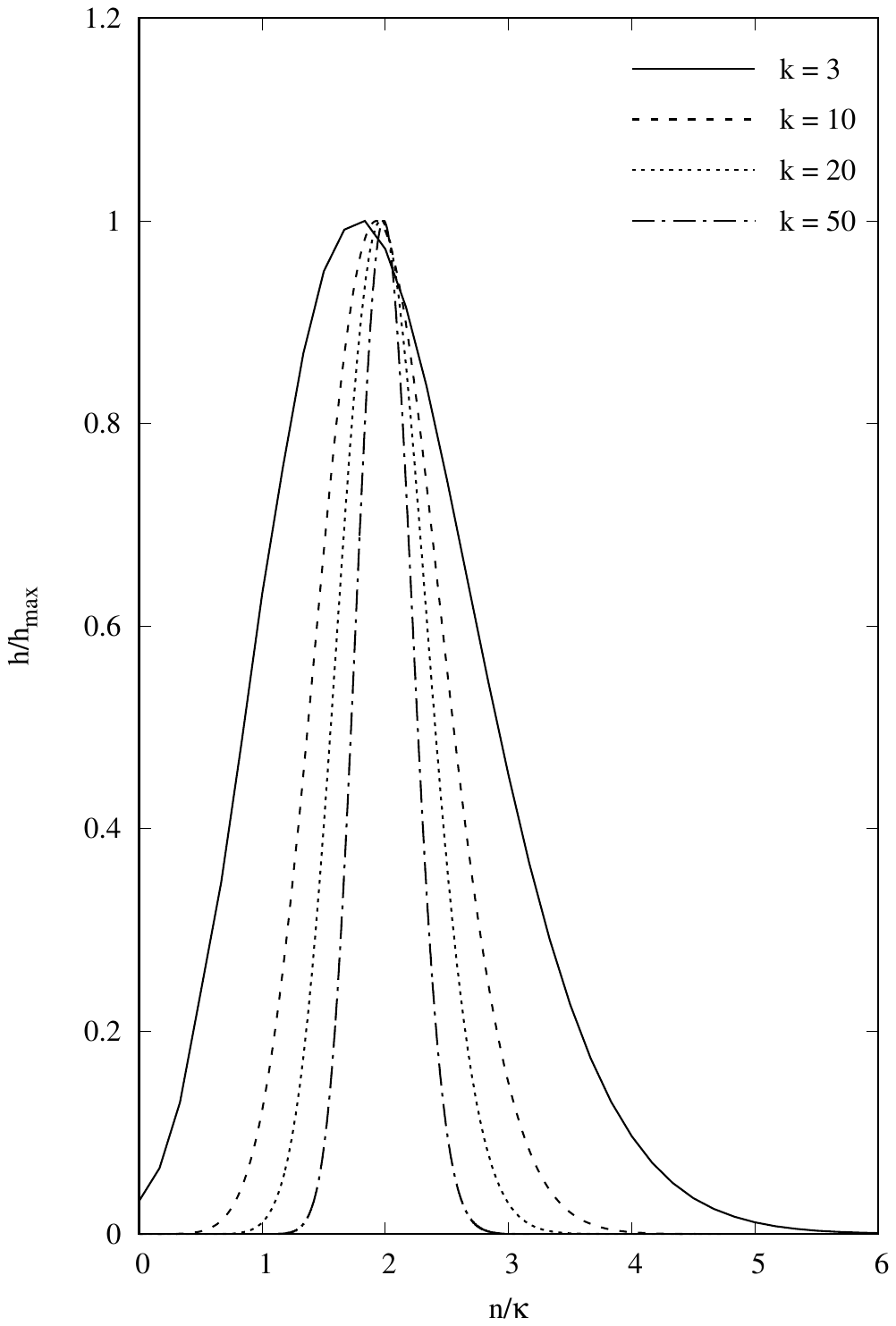}
\caption{\small
  \label{fig:hist_lam2}
  Plot of the pmf for the Poisson distribution of order $k$ for $\lambda=2$ and $k=3,10,20,50$.
  For each value of $k$, the pmf was scaled to a peak height of $1$.}
\end{figure}

\newpage
It was stated in \cite{Mane_Poisson_k_CC23_3} that the asymptotic formula for the mode is applicable for values $\lambda \ge 2$ (see eq.~\eqref{eq:Poisson_k_mode_int_n}).
Just out of curiosity, let us plot the pmf for $\lambda=2$ for the values $k=3,10,20,50$ employed in Sec.~\ref{sec:pmf}.
To plot them on the same scale, we scale all the curves to a peak height of $1$ and on the horizontal axis we plot the value of $n/\kappa \in [0,6]$.
The curves are displayed in Fig.~\ref{fig:hist_lam2}.
The mean is $\kappa\lambda=2\kappa$, hence the mean is at $n/\kappa=2$ for all the curves.
Observe that all the curves are smooth.
Issues of a local maximum at $n=k$, or a left peak, have faded out of the picture.
The smoothness improves as $k$ increases, which is expected.
The peaks get narrower as $k$ increases. We can explain this as follows.
Recall the variance is $\sigma_k^2(\lambda) = \frac16 k(k+1)(2k+1)\lambda$, hence in Fig.~\ref{fig:hist_lam2} the scaled peak width is (scaled standard deviation)
\bq
\begin{split}
\frac{\sigma_k(\lambda)}{\kappa} &= \sqrt{\frac{\frac16 k(k+1)(2k+1)\lambda}{\frac14 k^2(k+1)^2}}
\\
&= \sqrt{\frac{2\lambda}{3}}\, \sqrt{\frac{2k+1}{k(k+1)}} 
\\
&= \sqrt{\frac{2\lambda}{3}}\, \sqrt{\frac{1}{k}+\frac{1}{k+1}} \;.
\end{split}
\eq
Hence the scaled peak width decreases as the value of $k$ increases.
We fixed $\lambda=2$ in Fig.~\ref{fig:hist_lam2}.
Observe also that the peaks shift to the right as $k$ increases.
This is consistent with eq.~\eqref{eq:Poisson_k_mode_int_n}.
Set $n=\kappa\lambda=\mu_k(\lambda)$ and divide by $\kappa$ to obtain
\bq
\begin{split}
  \frac{\mu_k(\lambda) - m_k(\lambda)}{\kappa} &= \frac{\lfloor(3k+5)/8\rfloor}{\frac12 k(k+1)}
  \\
  &\le \frac{(3k+5)/8}{\frac12 k(k+1)}
  \\
  &= \frac14\, \biggl(\frac{5}{k} - \frac{2}{k+1}\biggr) \,.
\end{split}
\eq
The scaled difference between the mean and the mode decreases as the value of $k$ increases.

\newpage
\section{\label{sec:conc}Conclusion}
The major goal of this note was to characterize the structure of the probability mass function (pmf) of the Poisson distribution of order $k$.
The pmf can be partitioned into a single point at $n=0$, an increasing sequence for $n \in [1,k]$ and a mountain range for $n>k$.
That structure is by no means obvious from the formal definition of the pmf as a sum over terms with factorial denominators.
The mode structure of the pmf was quantified.
The locations of {\em and the reasons for} the mode jumps (where the mode increases by more than one unit) were established.
The ``excluded values'' tabulated in \cite{Mane_Poisson_k_CC23_5}, i.e.~integers which cannot be modes of the Poisson distribution of order $k$,
were explained as those integers which are skipped in the mode jumps.
It was demonstrated that for all $k\ge2$, the Poisson distribution of order $k$ has a denumerably infinite set of double modes, consisting of pairs of consecutive integers.
Numerical evidence was presented that the Poisson distribution of order $k$ does not have three or more joint modes.
The opportunity was also taken to publish improvements to various inequalities (sharper bounds, etc.)
and also to present new conjectured upper and lower bounds for the median and the mode of the Poisson distribution of order $k$.

%\newpage
%\section*{\label{sec:ack}Acknowledgements}

\newpage

\newpage

\newpage

\newpage

\newpage

\newpage

\newpage

\newpage

\newpage
\

\newpage

\end{document}